\documentclass{elsart3-1}

\usepackage{amssymb}

\usepackage[english,francais]{babel}
\usepackage{isolatin1}

\newtheorem{e-proposition}[theorem]{Proposition}

\newtheorem{e-definition}[theorem]{Definition\rm}

\newtheorem{theoreme}{Th\'eor\`eme}[section]
\newtheorem{lemme}[theoreme]{Lemme}
\newtheorem{proposition}[theoreme]{Proposition}

\newtheorem{definition}[theoreme]{D\'efinition\rm}

\renewcommand{\theequation}{\arabic{equation}}
\def\og{\leavevmode\raise.3ex\hbox{$\scriptscriptstyle\langle\!\langle$~}}
\def\fg{\leavevmode\raise.3ex\hbox{~$\!\scriptscriptstyle\,\rangle\!\rangle$}}

\def\<{\langle\,}
\def\>{\,\rangle}

\def\SG{{\mathfrak S}}

\def\MQSym{{\bf MQSym}}
\def\WQSym{{\bf WQSym}}
\def\PQSym{{\bf PQSym}}
\def\NCQSym{{\bf NCQSym}}

\def\tass{{\rm tass\,}}
\def\pack{{\rm tass\,}}
\def\detass{{\rm detass}}
\def\F{{\bf F}}
\def\G{{\bf G}}

\def\sl{\hbox{\goth sl\hskip 1pt}}

\def\K{{\mathbb K}}

\def\MM{{\mathcal M}} 
\def\M{{\bf M}}     

\def\a{{\bf a}}

\def\TT{{\mathcal T}}

\def\TD{{\mathfrak T}}
\def\limproj{{\rm lim\,proj\,}}

\begin{document}

\begin{frontmatter}




%
\selectlanguage{francais}
\title{Construction de trigèbres dendriformes}

\vspace{-2.6cm}
\selectlanguage{english}
\title{Construction of dendriform trialgebras}



\author[mlv]{Jean-Christophe Novelli},
\ead{novelli@univ-mlv.fr}
\author[mlv]{Jean-Yves Thibon}
\ead{jyt@univ-mlv.fr}
\address[mlv]{Institut Gaspard Monge, Universit\'e de
Marne-la-Vall\'ee, 77454 Marne-la-Vall\'ee Cedex 2, France}

\begin{abstract}
We realize the free dendriform trialgebra on one generator,
as well as several other examples of dendriform trialgebras,
as sub-trialgebras of an algebra of noncommutative polynomials
in infinitely many variables.
{\it To cite this article: J.-C. Novelli, J.-Y. Thibon, C. R.
Acad. Sci. Paris, Ser. I ??? (???).}

\vskip 0.5\baselineskip

\selectlanguage{francais}
\noindent{\bf R\'esum\'e}
\vskip 0.5\baselineskip
\noindent
Nous réalisons la trigèbre dendriforme libre sur un générateur, et
plusieurs autres exemples de trigèbres dendriformes, comme
sous-trigèbres d'une algèbre de polynômes non commutatifs en
une infinité de variables.
{\it Pour citer cet article~: J.-C. Novelli, J.-Y. Thibon,  C. R.
Acad. Sci. Paris, Ser. I ??? (???).}
\end{abstract}
\end{frontmatter}

\selectlanguage{english}
\section*{Abridged English version}
A {\em dendriform trialgebra}, as defined by
Loday and Ronco \cite{LR2}, is an associative algebra whose product splits
into three binary operations
$$
x\cdot y = x\prec y + x\circ y + x\succ y\,,
$$
where $\circ$ is associative, satisfying the compatibility conditions
(\ref{td1}) and (\ref{td2}).

The aim of this Note is to give explicit realizations of
various examples in terms of noncommutative polynomials in infinitely
many variables. Let $\K$ be a field of characteristic 0, and
$A=\{a_1<a_2 <\ldots\}$ be an infinite totally ordered alphabet. 
We denote by $A^*$ the free monoid over $A$ et $\K\< A \>=\limproj
\K\<A_n\>$ where $A_n$ is the interval $[a_1,a_n]$ of $A$ and $\K\<A_n\>$
the free associative algebra over $A_n$. 
We shall use the notations of \cite{NCSF6,NT2}.
We denote by $\max(w)$ the greatest
letter occuring in a word $w\in A^*$.

Our first result (Lemme \ref{KAtrid}) states that the ideal $\K\<A\>^+$
of polynomials without constant term is a tridendriform algebra for the
operations defined by Equations (\ref{ll}) to (\ref{rr}).

Next, we prove (Théorème \ref{TD}) that the sub-trialgebra $\TD$ of $\K\<A\>^+$
generated by the sum of the variables (\ref{M1}) is free as a dendriform
trialgebra.

This provides an explicit realization of this algebra, whose bases 
are known to be parametrized by {\em plane trees} (counted by
the little Schr\"oder numbers). To each such tree $T$, we associate
a polynomial $\MM_T$ (\ref{MT}). From each word
$w$ of length $n$, we build a plane tree $\TT(w)$ recursively defined
as follows. If $m=\max(w)$ and $w$ kas exactly $k$ occurences of $m$,
$\TT(w)$ is obtained from the factorization (\ref{fact}) by grafting
$\TT(v_0),\TT(v_1),\ldots,\TT(v_k)$ (in this order) on a common root.

The {\em packed word} $u=\pack(w)$ associated to a word $w\in A^*$ is
obtained by the following process. If
$b_1<b_2<\ldots <b_r$ are the letters occuring in $w$, $u$ is the
image of $w$ by the semigroup homomorphism $b_i\mapsto a_i$.
A word $u$ is said to be {\em packed} if $\pack(u)=u$.
To such a word, we associate a polynomial $\M_u$ (\ref{Mu})

These polynomials span a subalgebra of $\K\<A\>$, which consists in the
invariants of the noncommutative version of Hivert's quasi-symmetrizing
action (two words are in the same $\SG(A)$-orbit iff they have the
same packed word). Under the abelianization $\K\<A\>\rightarrow\K[X]$,
the $\M_u$ are mapped to the monomial quasi-symmetric functions $M_I$
($I$ being the evaluation vector of $u$). The algebra spanned
by the $\M_u$ is denoted by $\NCQSym$. It is known that this is a Hopf
subalgebra of $\MQSym$ \cite{HivT,NCSF6}.

Clearly, $\TD$ is contained in $\NCQSym$ (\ref{Mtu}). Moreover,
$\NCQSym$ can itself be embedded in the (dual) Hopf algebra
of parking functions $\PQSym^*$ of \cite{NT2} (Proposition \ref{wqsymG}), 
and its dual is explicitely embedded in $\PQSym$ (Proposition \ref{detass}).
The dual basis of $\M_u$ can be identified with $\F_\a$, where
$\a=\detass(u)$ is the {\em maximal unpacking} of $u$, that is, the
greatest (for the lexicographic order) parking function such that
$\tass(\a)=u$.

The next series of results makes use of Foissy's theory of 
{\em bidendriform bialgebras}  \cite{Foi}. We show that
$\NCQSym$, $\TD$, 
$\PQSym$ and $\MQSym$  are all bidendriform bialgebras, as well as
dendriform trialgebras (Théorèmes \ref{nctri}, \ref{tdbid},
\ref{pqbid} and \ref{mqbid}). 
As a consequence, all these are free and self-dual Hopf
algebras, and their primitive Lie algebras are free.

Finally, we observe that each homogeneous component $\NCQSym_n^*$
is stable under the internal product of parking functions introduced
in \cite{NT2}, and that the resulting algebra is isomorphic to
the Solomon-Tits algebra \cite{Tits} (Théorème \ref{solti}).

Théoreme \ref{sylv} provides an alternative construction of
the free dendriform trialgebra on one generator: it is isomorphic
to the quotient of $\NCQSym$ by the sylvester congruence of \cite{HNT}.

\selectlanguage{francais}
\section{Introduction}
\label{Intro}

Il est fréquent que la théorie des opérades et la théorie des fonctions
symétriques non commutatives conduisent, par des voies très différentes,
à la découverte des mêmes algèbres de Hopf basées sur des structures
combinatoires. Dans l'approche initiée par Loday \cite{Lod}, un problème
central est de décomposer la multiplication de chacune de ces algèbres en
somme de plusieurs opérations binaires, de manière à la faire apparaître
comme l'algèbre libre sur un générateur pour une  certaine opérade.
L'autre approche \cite{NCSF6,HNT,NTT,NT1,NT2} vise \`a faire apparaître ces
algèbres comme des généralisations de l'algèbre de Hopf des fonctions
symétriques, en les réalisant en termes de polynômes en un système auxiliaire
de variables, commutatives ou non, ce qui a habituellement pour effet de rendre
transparents la structure de Hopf et les divers morphismes généralement
présents dans ce contexte.

L'objet de cette Note est de proposer une telle réalisation des
trigèbres dendriformes de Loday et Ronco \cite{LR2}, et d'en tirer
quelques conséquences.

Dans la suite, $\K$ désignera un corps de caractéristique zéro. Nous utiliserons
les notations de \cite{NCSF6,NT1}.

\section{La trigèbre dendriforme libre sur un générateur}

Une {\em trigèbre dendriforme} est une algèbre associative dont la
multiplication se scinde en trois opérations
\begin{equation}
x\cdot y = x\prec y + x\circ y + x\succ y\,,
\end{equation}
où $\circ$ est associative, et
\begin{equation}\label{td1}
(x\prec y)\prec z = x\prec (y\cdot z)\,,\ \
(x\succ y)\prec z = x\succ (y\prec z)\,,\ \
(x\cdot y)\succ z = x\succ (y\succ z)\,,\ \
\end{equation}
\begin{equation}\label{td2}
(x\succ y)\circ z = x\succ (y\circ z)\,,\ \ \
(x\prec y)\circ z = x\circ (y\succ z)\,,\ \ \
(x\circ y)\prec z = x\circ (y\prec z)\,.
\end{equation}

Soit $A=\{a_1<a_2\ldots a_n <\ldots \,\}$ un alphabet dénombrable totalement
ordonné. On note $A^*$ le monoïde libre sur $A$ et $\K\< A \>=\limproj
\K\<A_n\>$ où $A_n$ est l'intervalle $[a_1,a_n]$ de $A$ et $\K\<A_n\>$
l'algèbre associative libre sur $A_n$. On notera $\max(w)$ la plus grande
lettre apparaissant dans un mot $w\in A^*$.

\medskip
\begin{definition}\label{deflcr}
Pour deux mots non vides $u,v\in A^*$, on pose
\begin{eqnarray}
u\succ v=\cases{uv &\mbox{si $\max(u)<\max(v)$}\cr 0 &\mbox{sinon,}}\label{ll}\\
u\circ v=\cases{uv &\mbox{si $\max(u)=\max(v)$}\cr 0 &\mbox{sinon,}}\label{cc}\\
u\prec v=\cases{uv &\mbox{si $\max(u)>\max(v)$}\cr 0 &\mbox{sinon.}}\label{rr}
\end{eqnarray}
\end{definition}

\begin{lemme}\label{KAtrid}
Les trois opérations $\prec,\circ,\succ$, munissent l'idéal d'augmentation $\K\<A\>^+$ d'une structure
de trigèbre dendriforme.
\end{lemme}

La vérification est immédiate. Considérons maintenant le polynôme
\begin{equation}\label{M1}
\M_1=\sum_{i\ge 1}a_i \,.
\end{equation}
Nous pouvons alors donner une réalisation explicite de la trigèbre
dendriforme
libre sur un générateur:

\medskip
\begin{theoreme}\label{TD}
La sous-trigèbre $\TD$ de  $\K\<A\>^+$ engendrée par $\M_1$ est libre en tant que trigèbre
dendriforme.
\end{theoreme}

\medskip
Pour établir ce résultat, il suffit de montrer que la série de Hilbert de $\TD$
(pour la graduation héritée de $\K\<A\>$) coïncide avec celle obtenue dans
\cite{LR2} pour la structure libre. Pour ce faire, associons à tout mot
$w$ de longueur $n$ un arbre plan $\TT(w)$ à $n+1$ feuilles, défini
récursivement comme suit. Si $m=\max(w)$ et si $w$ possède exactement $k$
occurences de $m$, écrivons
\begin{equation}\label{fact}
w=v_0\,m\,v_1\,m\,v_2\cdots v_{k-1}\,m\,v_k\,,
\end{equation}
où les $v_i$ sont éventuellement vides. Alors, $\TT(w)$ est l'arbre
formé des sous-arbres $\TT(v_0),\TT(v_1),\ldots,\TT(v_k)$ greffés (dans
cet ordre) sur une racine commune, avec la condition initiale
$\TT(\epsilon)=\emptyset$ pour le mot vide. On vérifie alors que
les polynômes 
\begin{equation}\label{MT}
\MM_T=\sum_{\TT(w)=T}w\,,
\end{equation}
où $T$ parcourt les arbres plans, sont dans la sous-trigèbre engendrée par
$\M_1$. \'Etant sur des ensembles de mots disjoints, ils sont linéairement
indépendants et le théorème s'ensuit. 

\medskip
On en déduit par exemple que la trigèbre dendriforme commutative libre sur un
générateur est isomorphe à $QSym^+$ \cite{LR2}.
En effet, c'est l'image de $\TD$ par le morphisme $\K\<A\>\rightarrow \K[X]$
qui envoie les $a_i$ sur des variables commutatives $x_i$.

\section{Partitions ordonnées et fonctions de parking}

Appelons {\em tassé} d'un mot $w$ le mot $\tass(w)$ obtenu par le procédé
suivant.
Soient $b_1\!<\!\ldots\!<\!b_i$ les lettres apparaissant dans $w$.
Alors, $\tass(w)$ est l'image de $w$ par le morphisme
$b_j\mapsto a_j$.

Un mot $u$ sera dit tassé si $\tass(u)=u$. Pour un tel mot, posons
\begin{equation}\label{Mu}
\M_u=\sum_{\tass(w)=u}w.
\end{equation}
On a clairement
\begin{equation}\label{Mtu}
\MM_T=\sum_{\TT(u)=T}\M_u\,.
\end{equation} 
Les mots tassés s'identifient naturellement aux partitions
ensemblistes ordonnées, et donc aux faces du permutoèdre, et
les arbres plans aux faces de l'associaèdre. Géométriquement, ces regroupements
correpondent à l'application de Tonks (\cite{Tonks}, Prop. 2.1).
De plus, il est connu que les $\M_u$ forment une base  
d'une sous-algèbre de $\K\<A\>$,
qui se plonge naturellement dans l'algèbre $\MQSym$ de \cite{NCSF6}, et hérite
ainsi d'une structure d'algèbre de Hopf. Les $\M_u$ sont les sommes des orbites
des mots  par l'action quasi-symétrisante d'Hivert \cite{HivT,Hiv}, et 
elles se projettent sur les fonctions quasi-symétriques monomiales
$M_I$
lorsqu'on fait commuter les variables. 
Nous la noterons $\NCQSym$ comme dans \cite{BZ} (elle était précédemment
notée $\WQSym$) o\`u les auteurs d\'emontrent qu'elle est libre et
colibre.

\medskip
\begin{theoreme}\label{nctri}
$\NCQSym$ est une trigèbre dendriforme, isomorphe à celle des partitions
ordonnées $\K[\Pi_\infty]$ de {\rm \cite{LR2}}.
De plus, c'est une bigèbre bidendriforme au sens de Foissy {\rm \cite{Foi}}.
\end{theoreme}

\medskip
\begin{theoreme}\label{tdbid}
$\TD$ est une sous-algèbre de Hopf de $\NCQSym$, également bidendriforme.
\end{theoreme}

\medskip
Il s'ensuit que $\NCQSym$ et $\TD$ sont autoduales, libres, colibres
et que leurs algèbres de Lie primitives sont libres.
On obtient ces résultats en plongeant
$\NCQSym$ dans l'algèbre de Hopf des fonctions de parking de \cite{NT1}.

\medskip
\begin{proposition}\label{wqsymG}
L'application linéaire $\NCQSym\rightarrow \PQSym^*$:
$\displaystyle \M_u\longmapsto \sum_{\tass(\a)=u}\G_\a$
est un monomorphisme d'algèbres de Hopf.
\end{proposition}

\medskip
\begin{theoreme}\label{pqbid}
$\PQSym$ est une bigèbre bidendriforme.
\end{theoreme}

\medskip
En particulier, $\PQSym$ est autoduale. 

Appelons {\em détassé maximal} $\a=\detass(u)$ d'un mot tassé $u$
la plus grande fonction de parking $\a$ (pour l'ordre lexicographique)
telle que $\tass(\a)=u$.
Un plongement naturel de $\NCQSym^*$ dans $\PQSym$ est donné par

\medskip
\begin{proposition}\label{detass}
Les $\F_\a$, où $\a$ parcourt les fonctions de parking détassées maximales,
engendrent une sous-algèbre de Hopf de $\PQSym$, isomorphe à $\NCQSym^*$. 
\end{proposition}

\medskip
Il est intéressant de noter que cette sous-algèbre est stable par le produit
intérieur défini dans \cite{NT2}. Rappelons que ce produit, noté $*$, préserve
chaque composante homogène $\PQSym_n$. 
La composante homogène $\NCQSym_n$ a même dimension que l'algèbre de
Solomon-Tits, définie dans \cite{Tits}, qui a pour support le complexe de
Coxeter de type $A_{n-1}$.
Le produit est  donné dans les deux cas par une formule explicite,
et par comparaison directe, on obtient

\medskip
\begin{theoreme}\label{solti}
La composante homogène $\NCQSym^*_n$, munie du produit intérieur au moyen de
la réalisation précédente, est isomorphe à l'algèbre de Solomon-Tits.
\end{theoreme}

\medskip
Il a été montré dans \cite{NT1} que le quotient de $\PQSym^*$ par la congruence
hypoplaxique (voir \cite{NCSF6}) était une algèbre de Hopf de même série de Hilbert
que $\TD$. Elle lui est isomorphe en tant qu'algèbre associative (car elle est
libre) mais pas en tant que cogèbre. En effet, sa duale contient une 
sous-algèbre
commutative isomorphe à $QSym$. Elle ne peut donc pas être autoduale.
Mais on peut tout de même réaliser $\TD$ par un procédé analogue :

\medskip
\begin{theoreme}\label{sylv}
Le quotient de l'algèbre $\NCQSym$, réalisée comme dans la Proposition \ref{wqsymG},
par la congruence sylvestre  {\rm \cite{HNT}} sur $A^*$, est isomorphe à
$\TD$.
\end{theoreme}

\medskip
En effet, on peut montrer que ce quotient hérite des structures de trigèbre
dendriforme et de bigèbre bidendriforme. 

Pour conclure, mentionnons également le résultat suivant, qui fournit une
nouvelle explication de l'autodualité de $\MQSym$ et de nombreuses autres
propriétés:

\medskip
\begin{theoreme}\label{mqbid}
$\MQSym$ est une bigèbre bidendriforme, et une trigèbre dendriforme.
\end{theoreme}



\begin{thebibliography}{00}
%
\bibitem{BZ} N. Bergeron, M. Zabrocki,
The Hopf algebras of non-commutative symmetric functions and quasi-symmetric
functions are free and cofree, math.CO/0509265.
%
\bibitem{NCSF6} G. Duchamp, F. Hivert, J.-Y. Thibon,
Noncommutative symmetric functions VI: free quasi-symmetric functions and
related algebras,
Internat. J. Alg. Comput. 12 (2002), 671--717.  
%
\bibitem{Foi} L. Foissy,
Bidendriform bialgebras, trees, and free quasi-symmetric functions,
math.RA/0505207.
%
\bibitem{HivT} F. Hivert, Combinatoire des fonctions quasi-symétriques,
Thèse de doctorat, Université de Marne-la-Vallée, 1999.
%
\bibitem{Hiv} F. Hivert, Hecke algebras, difference operators
and quasi-symmetric functions,
Advances in Math. 155 (2000), 181--238.
%
\bibitem{HNT} F. Hivert, J.-C. Novelli, J.-Y. Thibon,
The algebra of binary search trees, 
Theoret. Computer Sci. 339 (2005), 129--165.
%
\bibitem{Lod} J.-L. Loday, Scindement d'associativit\'e et alg\`ebres de
{H}opf,
Actes des Journ\'ees Math\'ematiques \`a la M\'emoire de Jean Leray,
S\'emin. Congr. Soc. Math. France 9 (2004), 155-172.
%
\bibitem{LR1} J.-L. Loday, M. O. Ronco,
Hopf algebra of the planar binary trees,
Adv. Math. 139 (1998), no. 2, 293--309.
%
\bibitem{LR2} J.-L. Loday, M. O. Ronco,
Trialgebras and families of polytopes,
Contemp. Math. {346} (2004), 369--398.
%
\bibitem{NT1} J.-C. Novelli, J.-Y. Thibon,
A Hopf algebra of parking functions,
FPSAC'04 (Vancouver, Juin 2004),
math.CO/0312126.        
%
\bibitem{NT2}  J.-C. Novelli, J.-Y. Thibon,
Parking functions and descent algebras, math.CO/0411387.
%
\bibitem{NTT} J.-C. Novelli, J.-Y. Thibon, N. M. Thiéry,
Algèbres de Hopf de graphes,   C. R. Math. Acad. Sci. Paris  339  (2004),  
607--610.
%
\bibitem{Tits} J. Tits,
Two properties of Coxeter complexes, Appendix to "A Mackey
formula in the group ring of a Coxeter group" (J. Algebra 41 (1976),
255--264) by Louis Solomon.  J. Algebra  41  (1976), 265--268. 
%
\bibitem{Tonks} A. Tonks, Relating the associahedron and the
permutohedron, 
Contemporary Math. 202 (1997), 33--36.
%
\end{thebibliography}
\end{document}